\documentclass{article}

\usepackage{amsmath}
\usepackage{amsthm}
\usepackage{alltt}

\usepackage[pdftex]{hyperref}
\usepackage{graphicx}

\newtheorem{theorem}{Theorem}[section]

\numberwithin{equation}{section}

\title{Toom-Cook Multiplication:\\ Some Theoretical and Practical Aspects}
\author{M.J. Kronenburg}

\begin{document}

\date{}
\maketitle

\begin{abstract}
Toom-Cook multiprecision multiplication is a well-known multiprecision
multiplication method,
which can make use of multiprocessor systems.
In this paper the Toom-Cook complexity is derived,
some explicit proofs of the Toom-Cook interpolation method are given,
the even-odd method for interpolation is explained,
and certain aspects of a 32-bit C++ and assembler implementation,
which is in development, are discussed.
A performance graph of this implementation is provided.
The Toom-Cook method can also be used to multithread other types
of multiplication, which is demonstrated for 32-bit GMP FFT multiplication.
\end{abstract}


\section{Integers as Polynomials}

Let the following polynomials be given:
\begin{equation}
 U(x) = \sum_{k=0}^n u_k x^k
\end{equation}
\begin{equation}
 V(x) = \sum_{k=0}^n v_k x^k
\end{equation}
\begin{equation}
 W(x) = \sum_{k=0}^{2n} w_k x^k
\end{equation}
Let $u$ and $v$ be two integers such that $u=U(2^b)$ and $v=V(2^b)$,
and let $w=uv$ be their product, so that $w=W(2^b)$.
This means that $u$ and $v$ are splitted into binary blocks $u_k$ and $v_k$ of $b$ bits wide.
Instead of multiplying $U(2^b)$ and $V(2^b)$ directly to obtain $W(2^b)$,
Toom-Cook multiplication computes $y_k=W(x_k)=U(x_k)V(x_k)$ for $2n+1$ distinct values of $x_k$,
obtains the coefficients $w_k$ of $W(x)$ by polynomial interpolation at the points $y_k=W(x_k)$,
and obtains the result $uv=W(2^b)$.
The values of $x_k$ are very small and chosen in such a way that $U(x_k)$ and $V(x_k)$
can be computed efficiently and interpolation is efficient.
The products $U(x_k)V(x_k)$ are independent,
and can therefore be computed in parallel on multiple processors,
using Toom-Cook or any other method of multiplication.

\section{Toom-Cook Complexity}
From the definition of Toom-Cook multiplication above its complexity can be derived,
assuming that the products are also recursively computed with Toom-Cook multiplication,
and assuming no parallel processing.
Let the time to perform the Toom-Cook multiplication be order $T(N)$,
where $N$ is the size of the integers (in bits, bytes, words or decimals),
and let $B=n+1$ be the number of $b$-bit blocks of $u$ and $v$,
which is called $B$-way Toom-Cook multiplication.
Then the number of multiplications that is needed in the algorithm is $M=2n+1=2B-1$.
Therefore the following relation exists:
\begin{equation}
 T(N) = M\cdot T(\frac{N}{B})
\end{equation}
By assuming that $T(N)=N^{\displaystyle \alpha}$:
\begin{equation}
 N^{\displaystyle \alpha} = M \left(\frac{N}{B}\right)^{\displaystyle \alpha}
\end{equation}
By taking logarithms on both sides:
\begin{equation}
 \alpha = \frac{\log(M)}{\log(B)} = \frac{\log(2B-1)}{\log(B)}
\end{equation}
and the conclusion is that the complexity of $B$-way Toom-Cook is order:
\begin{equation}
 T(N) = N^{\displaystyle\frac{\log(2B-1)}{\log(B)}}
\end{equation}
For Karatsuba multiplication, which is 2-way Toom-Cook,
$B=2$ and $T(N)=N^{1.585}$.
For 16-way Toom-Cook this would be $T(N)=N^{1.239}$,
and for 32-way Toom-Cook this would be $T(N)=N^{1.195}$.
This complexity is worse than FFT (fast Fourier transform)
multiplication,
but the strong point of Toom-Cook, that is parallel processing,
has not yet been taken into account.
When there are $P$ parallel processors which are only used for
the top-level products so that $P\leq 2B-1$:
\begin{equation}
 T(N) = \frac{2B-1}{P}\left(\frac{N}{B}\right)^{\displaystyle\frac{\log(2B-1)}{\log(B)}}
 = \frac{1}{P}\cdot N^{\displaystyle\frac{\log(2B-1)}{\log(B)}}
\end{equation}
When enough parallel processors and memory are available,
the lower-level products can also be done in parallel, making the algorithm even faster.

\newpage
\section{Polynomial Interpolation}

The interpolation problem is how to find the coefficients $c_l$, $0\leq l\leq n$,
of a polynomial of degree $n$, $p_n(x)=\sum_{l=0}^nc_lx^l$,
given $n+1$ distinct points $p_n(x_k)=y_k$ for $0\leq k\leq n$.
The following theorem is given in exercise 15 of chapter 4.6.4 of \cite{K97}:
\begin{theorem}
Let $p_n(x)$ be a polynomial in $x$ of degree $n$,
and let $n+1$ distinct values of this polynomial be given by $p_n(x_k)=y_k$ for $0\leq k\leq n$.
Then this polynomial is given by Newton interpolation:
\begin{equation}\label{newtoninterp}
  p_n(x) = \sum_{k=0}^n \alpha_k \prod_{j=0}^{k-1} (x-x_j)
\end{equation}
where \cite{K97}:
\begin{equation}\label{newtonalpha}
 \alpha_k = \sum_{i=0}^k \frac{y_i}{\displaystyle\prod_{\substack{j=0\\j\neq i}}^k (x_i-x_j)}
\end{equation}
\end{theorem}
\begin{proof}
Substituting (\ref{newtonalpha}) in (\ref{newtoninterp}),
changing the order of summation and interchanging $k$ and $i$ yields:
\begin{equation}
 p_n(x) = \sum_{k=0}^n y_k \sum_{i=k}^n \frac{\displaystyle\prod_{j=0}^{i-1}(x-x_j)}
   {\displaystyle\prod_{\substack{j=0\\j\neq k}}^i(x_k-x_j)}
\end{equation}
The Lagrange-Waring interpolation \cite{K97} of this polynomial for these points is:
\begin{equation}
  p_n(x) = \sum_{k=0}^n y_k \prod_{\substack{j=0\\j\neq k}}^n \frac{x-x_j}{x_k-x_j}
\end{equation}
Equating each term of the last two identities for each $k$,\newline
multiplying both sides by $\displaystyle\prod_{\substack{j=0\\j\neq k}}^n(x_k-x_j)$
and dividing both sides by $\displaystyle\prod_{j=0}^{k-1}(x-x_j)$ yields:
\begin{equation}
  \prod_{j=k+1}^n (x-x_j) = \sum_{i=k}^n \prod_{j=k}^{i-1}(x-x_j) \prod_{j=i+1}^n(x_k-x_j)
\end{equation}
For each $k$, by renumbering the $x_j$, $k$ can be replaced by $0$ and $n$ by $n-k$,
and splitting off the $i=0$ term, and using $x-x_0=(x-x_i)-(x_0-x_i)$ yields:
\begin{equation}
\begin{split}
  \prod_{j=1}^n(x-x_j) & - \prod_{j=1}^n(x_0-x_j) = (x-x_0)\sum_{i=1}^n\prod_{j=1}^{i-1}(x-x_j)\prod_{j=i+1}^n(x_0-x_j) \\
 & = \sum_{i=1}^n\left[\prod_{j=1}^i(x-x_j)\prod_{j=i+1}^n(x_0-x_j) - \prod_{j=1}^{i-1}(x-x_j)\prod_{j=i}^n(x_0-x_j)\right]
\end{split}
\end{equation}
The last sum is a telescoping series:
on the right side the first product for each $i$ cancels the second product
for each $i+1$, only leaving the first product for $i=n$ and the second product for $i=1$.
These two remaining products are the two products on the left side.
\end{proof}
For computing the $\alpha_k$ the following theorem is used \cite{K97}:
\begin{theorem}
Let for $0\leq m\leq n-k$ the $\alpha_k^{(m)}$ be given by:
\begin{equation}
  \alpha_k^{(m)} = \sum_{i=0}^k \frac{y_{m+i}}{\displaystyle\prod_{\substack{j=0\\j\neq i}}^k (x_{m+i}-x_{m+j})}
\end{equation}
Then \cite{K97}:
\begin{equation}
  \alpha_k^{(m)} = \frac{\alpha_{k-1}^{(m+1)}-\alpha_{k-1}^{(m)}}{x_{m+k}-x_m}
\end{equation}
\end{theorem}
\begin{proof}
\begin{equation}
  \alpha_{k-1}^{(m+1)} = \sum_{i=0}^{k-1}\frac{y_{m+i+1}}{\displaystyle\prod_{\substack{j=0\\j\neq i}}^{k-1}(x_{m+i+1}-x_{m+j+1})}
   = \sum_{i=1}^{k}\frac{y_{m+i}}{\displaystyle\prod_{\substack{j=1\\j\neq i}}^{k}(x_{m+i}-x_{m+j})}
\end{equation}
\begin{equation}
  \alpha_{k-1}^{(m)} = \sum_{i=0}^{k-1}\frac{y_{m+i}}{\displaystyle\prod_{\substack{j=0\\j\neq i}}^{k-1}(x_{m+i}-x_{m+j})}
\end{equation}
\begin{equation}
\begin{split}
  \alpha_{k-1}^{(m+1)}-\alpha_{k-1}^{(m)} & = \sum_{i=1}^{k-1}\frac{y_{m+i}}{\displaystyle\prod_{\substack{j=1\\j\neq i}}^{k-1}(x_{m+i}-x_{m+j})}
   \left[\frac{1}{x_{m+i}-x_{m+k}}-\frac{1}{x_{m+i}-x_m}\right] \\
  & \quad + \frac{y_{m+k}}{\displaystyle\prod_{j=1}^{k-1}(x_{m+k}-x_{m+j})}
    - \frac{y_m}{\displaystyle\prod_{j=1}^{k-1}(x_m-x_{m+j})} 
\end{split}
\end{equation}
Using:
\begin{equation}
  \frac{1}{x_{m+i}-x_{m+k}}-\frac{1}{x_{m+i}-x_m} = \frac{x_{m+k}-x_m}{(x_{m+i}-x_{m+k})(x_{m+i}-x_m)}
\end{equation}
this becomes:
\begin{equation}
\begin{split}
  \alpha_{k-1}^{(m+1)}-\alpha_{k-1}^{(m)} 
  & = (x_{m+k}-x_m)\sum_{i=1}^{k-1}\frac{y_{m+i}}{\displaystyle\prod_{\substack{j=0\\j\neq i}}^k(x_{m+i}-x_{m+j})} \\
  & \quad + (x_{m+k}-x_m)\Bigg[\frac{y_{m+k}}{\displaystyle\prod_{j=0}^{k-1}(x_{m+k}-x_{m+j})}
                             + \frac{y_m}{\displaystyle\prod_{j=1}^k(x_m-x_{m+j})}\Bigg] \\
  & = (x_{m+k}-x_m) \sum_{i=0}^k\frac{y_{m+i}}{\displaystyle\prod_{\substack{j=0\\j\neq i}}^k(x_{m+i}-x_{m+j})} \\
  & = (x_{m+k}-x_m) \alpha_k^{(m)}
\end{split}
\end{equation}
\end{proof}
Using this formula, $\alpha_0^{(m)}=y_m$ and $\alpha_k^{(0)}=\alpha_k$.
When mapping the values of $m$ to an array, because $0\leq m\leq n-k$, results are overwritten.
To avoid this, in the loop, $k$ is added to every $m$, so that $k\leq m\leq n$,
and the loop direction is reversed, so that the resulting $\alpha_k$ remain in the array.
This results in the following algorithm \cite{K97}:
\begin{alltt}
for(k=0;k<=n;k++)
 coeff[k] = y[k];
for(k=1;k<=n;k++)
 for(m=n;m>=k;m--)
  coeff[m] = (coeff[m]-coeff[m-1])/(x[m]-x[m-k]);
\end{alltt}
where the divisions are exact integer divisions \cite{GM94}.\newline
For getting from the $\alpha_k$ to the coefficients,
let $u_n^{(m)}(x)$ be the following intermediate polynomial:
\begin{equation}
 u_n^{(m)}(x) = \sum_{k=0}^{m-1} \beta_k^{(m)}\prod_{j=0}^{k-1}(x-x_j)
   + \sum_{k=m}^n \beta_k^{(m)}x^{k-m}\prod_{j=0}^{m-1}(x-x_j)
\end{equation}
so that:
\begin{equation}
 u_n^{(n)}(x) = \sum_{k=0}^n \beta_k^{(n)} \prod_{j=0}^{k-1}(x-x_j)
\end{equation}
and:
\begin{equation}
 u_n^{(0)}(x) = \sum_{k=0}^n \beta_k^{(0)} x^k
\end{equation}
This means that the $\beta_k^{(n)}$ are the $\alpha_k$ and the $\beta_k^{(0)}$ are
the coefficients.
\begin{theorem}
For getting from $u_n^{(m+1)}(x)$ to $u_n^{(m)}(x)$ the following relation can be used
for $k$ from $m$ to $n-1$:
\begin{equation}
 \beta_k^{(m)} = \beta_k^{(m+1)} - x_m \beta_{k+1}^{(m+1)}
\end{equation}
\end{theorem}
\begin{proof}
$u_n^{(m+1)}(x)$ is given by:
\begin{equation}
 u_n^{(m+1)}(x) = \sum_{k=0}^m \beta_k^{(m+1)}\prod_{j=0}^{k-1}(x-x_j)
   + \sum_{k=m+1}^n \beta_k^{(m+1)}x^{k-m-1}\prod_{j=0}^m(x-x_j)
\end{equation}
Then it is clear that $u_n^{(m+1)}(x)$ becomes $u_n^{(m)}(x)$
by repetitive use for $k$ from $m$ to $n-1$ of:
\begin{equation}
\begin{split}
 & \beta_k^{(m+1)}x^{k-m}\prod_{j=0}^{m-1}(x-x_j) + \beta_{k+1}^{(m+1)}x^{k-m}\prod_{j=0}^m(x-x_j) \\
 & = (\beta_k^{(m+1)} - x_m \beta_{k+1}^{(m+1)}) x^{k-m}\prod_{j=0}^{m-1}(x-x_j)
  + \beta_{k+1}^{(m+1)}x^{k-m+1}\prod_{j=0}^{m-1}(x-x_j) \\
 & = \beta_k^{(m)} x^{k-m}\prod_{j=0}^{m-1}(x-x_j)  + \beta_{k+1}^{(m+1)}x^{k-m+1}\prod_{j=0}^{m-1}(x-x_j)
\end{split}
\end{equation}
\end{proof}
This results in the following algorithm \cite{K97}:
\begin{alltt}
for(m=n-1;m>=0;m--)
 for(k=m;k<=n-1;k++)
  coeff[k] -= x[m] * coeff[k+1];
\end{alltt}
For obtaining the polynomial coefficients from the $x_k$ and $y_k$
these two double loops are executed, and because the $x_k$ are very small,
they have complexity order $n^2Y$ where $Y$ is the size of the $y_k$.
Because $y_k=U(x_k)V(x_k)$, and the $x_k^n$ are also small,
and the size of the binary blocks is about $b=N/B$ where $B=n+1$,
the size $Y$ is about $2N/(n+1)$.
Therefore the complexity of both double loops becomes order $nN$.
The complexity of computing the $U(x_k)$ and $V(x_k)$ is order $n^2b$, so this complexity is also order $nN$.
The complexity of the total overhead of Toom-Cook multiplication is therefore order $nN$.
Because the multiplication complexity is order $N^\alpha$ where $\alpha>1$ (see above), for constant $B$,
the percentage of time used for overhead decreases with increasing $N$.\\
Another approach is to solve the $c_l$, $0\leq l\leq n$, from $p_n(x_k)=y_k=\sum_{l=0}^n c_l x_k^l$ for $0\leq k\leq n$
by putting the $x_k^l$ in a square matrix and inverting this matrix \cite{BZ,Ggmpman}.
Then for specific $n$ the computation of the $c_l$ can be optimized from the entries of this inverse matrix.
For specific $n$ this may be faster, but for an algorithm that works for any $n$
the more general algorithm above may be preferred.

\section{The Even-Odd Method for Interpolation}
The even-odd method for interpolation is introduced in exercise 4 of chapter 4.3.3 of \cite{K97}.
The polynomials $U(x)$, $V(x)$ and $W(x)$ can be splitted into parts with
even and odd powers \cite{K97}:
\begin{equation}
 U(x) = \sum_{k=0}^n u_k x^k = \sum_{k=0}^{\lfloor n/2\rfloor} u_{2k} x^{2k}
   + x \sum_{k=0}^{\lfloor (n-1)/2\rfloor} u_{2k+1} x^{2k} 
 = U_e(x^2) + x U_o(x^2)
\end{equation}
and likewise $V(x)=V_e(x^2)+x V_o(x^2)$ and:
\begin{equation}\label{wsum}
 W(x) = \sum_{k=0}^{2n} w_k x^k = \sum_{k=0}^{n} w_{2k} x^{2k}
   + x \sum_{k=0}^{n-1} w_{2k+1} x^{2k} 
 = W_e(x^2) + x W_o(x^2)
\end{equation}
Then $W(x)$ and $W(-x)$ are calculated:
\begin{equation}\label{wpos}
  W(x) = [U_e(x^2)+x U_o(x^2)][V_e(x^2)+x V_o(x^2)]
\end{equation}
\begin{equation}\label{wneg}
  W(-x) = [U_e(x^2)-x U_o(x^2)][V_e(x^2)-x V_o(x^2)]
\end{equation}
and $W_e(x^2)$ and $W_o(x^2)$ are obtained:
\begin{equation}\label{wodd}
 W_e(x^2) = \frac{1}{2}(W(x)+W(-x))
\end{equation}
\begin{equation}\label{weven}
 W_o(x^2) = \frac{1}{2x}(W(x)-W(-x))
\end{equation}
Because the total number of multiplications $M=2B+1$ is uneven,
taking $x_0=0$, $W_e(0)=U(0)V(0)$ is calculated separately.
For using the interpolation algorithm with identical parameters for $W_e(x^2)$ and $W_o(x^2)$,
the $W_o(x^2)$ should be made of degree $n$ instead of $n-1$.
This is done by multiplying $W_o(x^2)$ by $x^2$ and adding as first point $x=0,y=0$,
so that its array positions are aligned to those of $W_e(x^2)$,
and the resulting coefficients move one place higher.
This way the two polynomials $W_e(x^2)$ and $x^2W_o(x^2)$ of degree $n$
can be interpolated instead of $W(x)$ of degree $2n$.
As the $x_k$ are very small, so are the $x_k^2$, and therefore
the complexity of interpolation is order $n^2Y$, where $Y$ is the size of the $y_k$,
see above. This size $Y$ is about identical with or without the even-odd method.
The extra complexity of (\ref{wpos}) to (\ref{weven}) is order $nY$ which can be neglected.
So without the even-odd method the interpolation time is $(2n)^2Y$
and with the even-odd method it is $2n^2Y$.
Therefore the even-odd method reduces the time of interpolation by a factor $2$.
For large multiplications, because the interpolations of these two polynomials
are independent, they can be done in parallel, reducing the time of interpolation
by another factor 2.

\newpage
\section{Using Toom-Cook for Multithreading GMP FFT Multiplication}

When the Toom-Cook algorithm is used for multithreading, then the
multiplications in each thread may be other Toom-Cook multiplications,
but may also be any other type of multiplications.
A good candidate is the FFT multiplication of GMP (version 5.0.2) \cite{Ggmp,Ggmpman}.
Assuming that the complexity of GMP FFT multiplication is order $N\log(N)$,
and when there are $P$ parallel processors so that $P\leq 2B-1$,
then the multiplication time $T(N)$ is about:
\begin{equation}
 T(N) = \frac{2B-1}{P}\frac{N}{B}\log(\frac{N}{B}) \simeq \frac{2B-1}{PB}N\log(N)
\end{equation}
When there are enough processors, so that $P=2B-1$, we get an improvement factor
of $1/B$, where for pure Toom-Cook this was $1/P$.\\
Instead of first choosing the $B$ for Toom-Cook multiplication and then dividing
the multiplications over the threads, now
the $B$ is chosen such that $2B-1$ fits into the number of threads,
that is $B$ is always half of the number of threads.
This way the amount of Toom-Cook overhead is also minimized.

\section{An Implementation}

\begin{figure}[h]
\caption{Performance of 32-bit 16-way Toom-Cook multiplication}
\begin{center}
\includegraphics[width=10.5cm,height=14cm,angle=-90,trim=0 1cm 0 0,keepaspectratio]{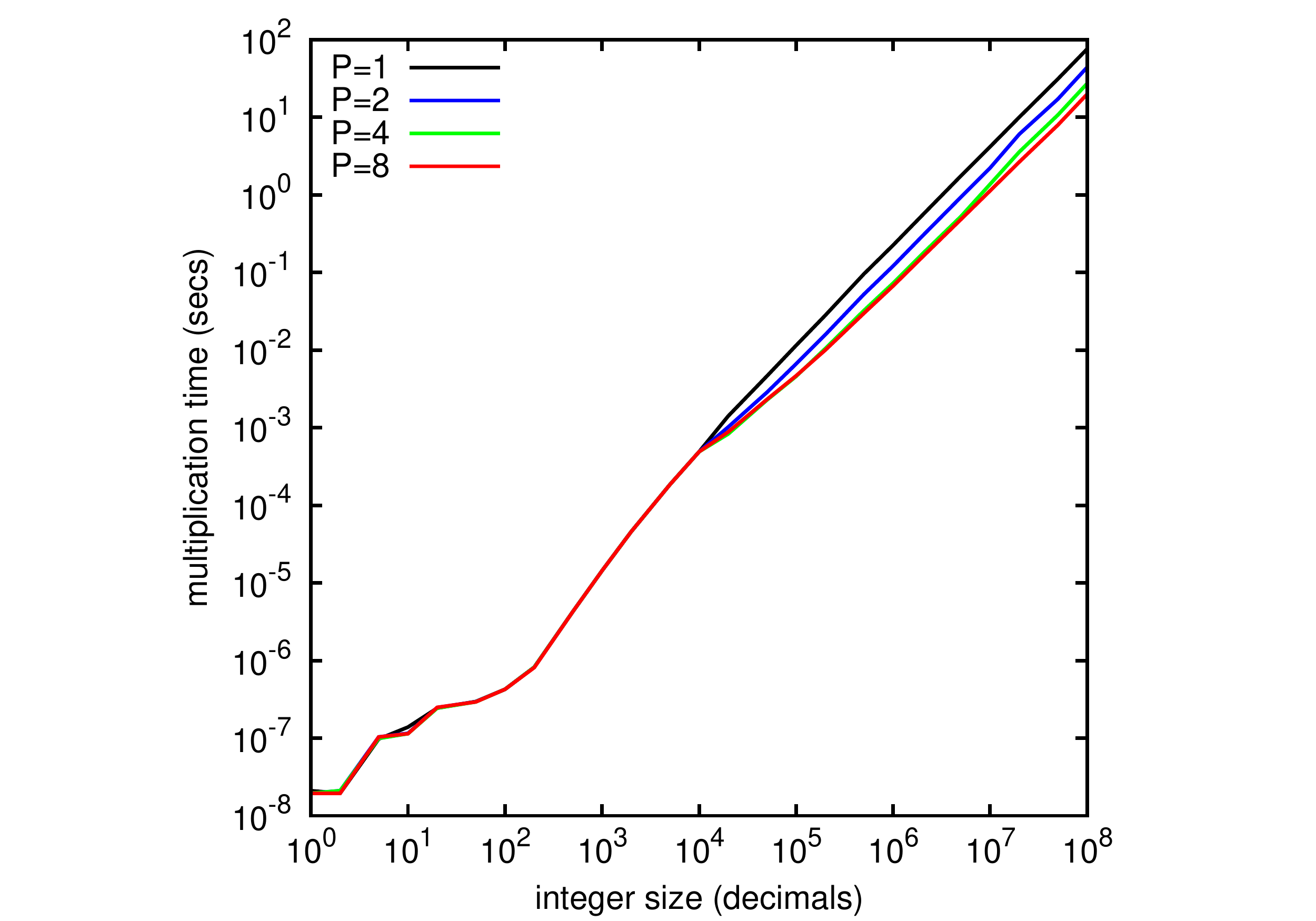}
\end{center}
\label{fig1}
\end{figure}

\begin{figure}[t]
\caption{Performance of 32-bit Toom-Cook with GMP FFT multiplication}
\begin{center}
\includegraphics[width=10.5cm,height=14cm,angle=-90,trim=0 1cm 0 0,keepaspectratio]{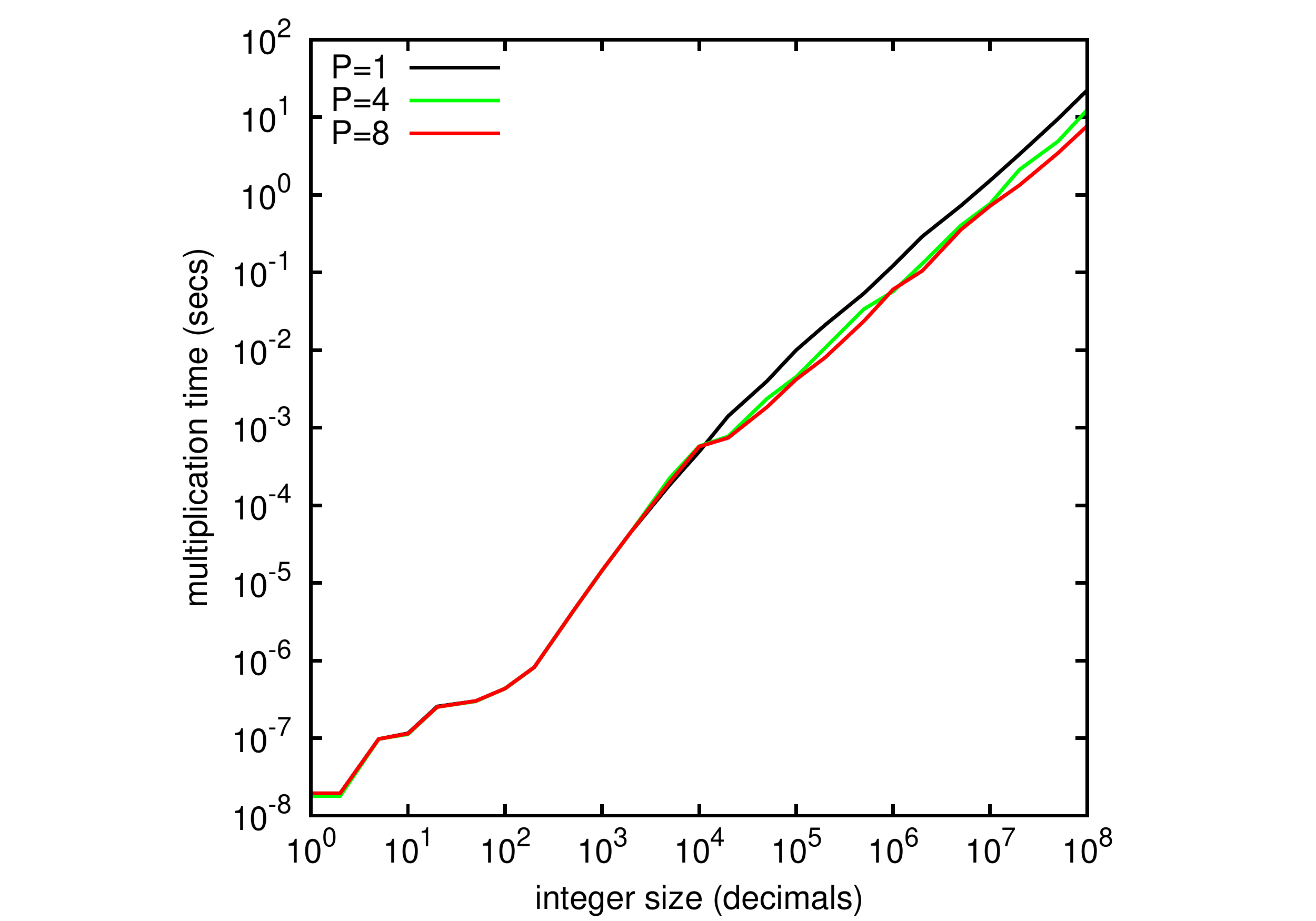}
\end{center}
\label{fig2}
\end{figure}

A 32-bit C++ and assembler implementation of Toom-Cook multiplication was developed.
The maximum parameter $B=16$ so that $n=15$, and the interpolation points $x_k$ were 
chosen as $0,1,2,4,...,2^{n-1}$ so that many multiplications reduce to shifts \cite{K97}.
This way the maximum $x_k^2$ was $2^{2(n-1)}=2^{28}$ which fits in a 32-bit word,
so that the exact divisions in the interpolations are all shifts or
exact divisions by a single $32$-bit word \cite{GM94}.
The program was compiled and run on an Intel Core i7 2.67 GHz machine which has 4 processors
and on which can run 8 threads.
The performance of its multiplication is shown in table \ref{tab2} and figure \ref{fig1}.
The multiplications in table \ref{tab2} were for $10^8$ decimals and the timings were in seconds.
The top-level overhead time is given and is included in the total time.
In figure \ref{fig1} the Toom-Cook multiplication starts at about 3000 decimals,
and its multithreading starts at about 13000 decimals.
For large multiplications the slope is 1.265,
which is close to the theoretical 1.239 (see above).\\
The extra memory required by Toom-Cook multiplication is $6$ times the size
of the argument when no multithreading is used, but when using $nthr$ top-level threads
this memory may increase with $6\cdot nthr/B$ times the size of the argument.
The size of an argument with $d$ decimals is about $d\cdot \ln(10)/(8\cdot\ln(2))$ bytes,
which for $d=10^8$ is about 42 Mbytes.\\
32-bit Toom-Cook with 32-bit GMP FFT multiplication \cite{Ggmp,Ggmpman} with $10^8$ decimals on 8 threads
is about 3 times faster than pure 32-bit GMP FFT multiplication on 1 thread,
see table \ref{tab2} and figure \ref{fig2}.\\
In table \ref{tab1} the full 8 threads were always used,
and bindec means conversion from binary to decimal when a result is written to output or to file.
The decimals of the constants were all checked with Mathematica$^{\textregistered}$ \cite{W03}.
The performance would improve by using 64-bit C++ and assembler
(which allows $B=32$) and using more processors.

\hypertarget{tabs}{}
\pdfbookmark[0]{Tables}{tabs}

\begin{table}[h]
\caption{Timings of 32-bit Toom-Cook with/without GMP FFT: multiplication}\label{tab2}
\begin{center}
\begin{tabular}{ | l | l | l | l | l | }
\hline
P & B & Toom-Cook & improvement & overhead \\ \hline
1 & 16 & 75.8 & 1.00 & 3.5 \\ \hline
2 & 16 & 44.1 & 1.72 & 2.4 \\ \hline
4 & 16 & 27.2 & 2.79 & 2.4 \\ \hline
8 & 16 & 19.9 & 3.81 & 2.3 \\ \hline
 & & Toom-Cook with GMP FFT & & \\ \hline
1 &    & 22.3 & 1.00 & \\ \hline
4 & 2  & 12.4 & 1.80 & 0.24 \\ \hline
8 & 4  & 7.7  & 2.90 & 0.53 \\ \hline
\end{tabular}
\end{center}
\end{table}

\begin{table}[h]
\caption{Timings of 32-bit Toom-Cook with/without GMP FFT: constants}\label{tab1}
\begin{center}
\begin{tabular}{ | l | l | l | l | l | l | l | }
\hline
 & \multicolumn{6}{c|}{Toom-Cook} \\ \hline
 decimals & $\sqrt{2}$ & $e$ & $\pi$ & $\log(2)$ & $\gamma$ & bindec \\ \hline 
$10^5$ & 0.026 & 0.059 & 0.14 & 0.27 & 2.84 & 0.063 \\ \hline
$10^6$ & 0.24  & 0.54 & 1.62 & 4.07 & 63.7 & 0.76   \\ \hline
$10^7$ & 2.83  & 9.00 & 30.5 & 80.4 & 1595 & 11.2   \\ \hline
$10^8$ & 46.1  & 157 & 594 & 1620  &      & 190    \\ \hline
 & \multicolumn{6}{c|}{Toom-Cook with GMP FFT} \\ \hline
$10^5$ & 0.023 & 0.055 & 0.14 & 0.27 & 2.70 & 0.059 \\ \hline
$10^6$ & 0.21  & 0.54 & 1.56 & 3.78 & 48.5 & 0.73   \\ \hline
$10^7$ & 2.02  & 6.67 & 22.2 & 58.9 & 950  & 9.14   \\ \hline
$10^8$ & 21.1  & 80.9 & 302  & 794  &      & 110    \\ \hline
\end{tabular}
\end{center}
\end{table}

\hypertarget{refs}{}
\pdfbookmark[0]{References}{refs}

\hypertarget{figs}{}
\pdfbookmark[0]{Figures}{figs}

\end{document}